\documentclass[12pt]{amsart}
\usepackage{amscd}
\usepackage{amssymb}

\newcommand{\nc}{\newcommand}
\newcommand{\rnc}{\renewcommand}
\nc{\nth}{\newtheorem}

\nth{thm}{Theorem}
\nth{pr}{Proposition}
\nth{lem}{Lemma}
\nth{cor}{Corollary}

\theoremstyle{definition}

\nth{rem}{Remark}

\nth{df}{Definition}

\nc{\CCC}{{\mathbb C}}
\nc{\PP}{{\mathbb P}}
\nc{\ZZ}{{\mathbb Z}}
\nc{\BA}{{\mathbb A}}
\nc{\BN}{{\mathbb N}}
\nc{\BZ}{{\mathbb Z}}

\nc{\bG}{{\mathbf G}}
\nc{\bL}{{\mathbf L}}
\nc{\bX}{{\mathbf X}}
\nc{\bB}{{\mathbf B}}
\nc{\bP}{{\mathbf P}}
\nc{\bH}{{\mathbf H}}

\nc{\CM}{{\mathcal M}}
\nc{\CO}{{\mathcal O}}

\nc{\fb}{{\mathfrak b}}
\nc{\fg}{{\mathfrak g}}
\nc{\fh}{{\mathfrak h}}
\nc{\fn}{{\mathfrak n}}
\nc{\fp}{{\mathfrak p}}
\nc{\fr}{{\mathfrak r}}
\nc{\fsl}{{\mathfrak sl}}

\nc{\al}{\alpha}
\nc{\la}{{\check\lambda}}
\rnc{\o}{{\check\omega}}
\nc{\Om}{\mho}

\nc{\refeq}[1]{$(\ref{#1})$}
\nc{\om}{\otimes\CO_C(-1)}
\nc{\opl}{\mathop{\oplus}\limits}

\rnc{\tilde}{\widetilde}
\rnc{\hat}{\widehat}

\begin{document}

\title[A note on the symplectic structure on the space of $G$-monopoles]
{A note on the symplectic structure\\ on the space of $G$-monopoles}

\author[M.~Finkelberg]{Michael Finkelberg}
\address{\!\!\!\!\!\!\!\!\! Independent Moscow University,
11 Bolshoj Vlasjevskij pereulok,
Moscow 121002 Russia}
\email{fnklberg@mccme.ru}
\author[A.~Kuznetsov]{Alexander Kuznetsov}
\address{\!\!\!\!\!\!\!\!\! Independent Moscow University,
11 Bolshoj Vlasjevskij pereulok,
Moscow 121002 Russia}
\email{sasha@kuznetsov.mccme.ru kuznetsov@mpim-bonn.mpg.de}
\author[N.~Markarian]{Nikita Markarian}
\address{\!\!\!\!\!\!\!\!\! Independent Moscow University,
11 Bolshoj Vlasjevskij pereulok,
Moscow 121002 Russia}
\email{nikita@mccme.ru}
\author[I.~Mirkovi\'c]{Ivan Mirkovi\'c}
\address{\!\!\!\!\!\!\!\!\! Dept. of Mathematics and Statistics,
University of Massachusetts at Amherst, Amherst MA 01003-4515, USA}
\email{mirkovic@math.umass.edu}
\thanks{M.F. and A.K. were partially supported by CRDF grant RM1-265}
\dedicatory{To Sasha Shen on his 40th birthday}
\begin{abstract}
Let $\bG$ be a semisimple complex Lie group with a Borel subgfoup $\bB$.
Let $\bX=\bG/\bB$ be the flag manifold of $\bG$. Let $C=\PP^1\ni\infty$ be
the projective line. Let $\alpha\in H_2(\bX,\BZ)$. The moduli space of
$\bG$-monopoles of topological charge $\alpha$ is naturally identified with
the space $\CM_b(\bX,\alpha)$ of based maps from $(C,\infty)$ to $(\bX,\bB)$
of degree $\alpha$. The moduli space of $\bG$-monopoles carries a natural
hyperk\"ahler structure, and hence a holomorphic symplectic structure. It was
recently explicitly computed by R.~Bielawski in case $\bG=SL_n$. We propose
a simple explicit formula for another natural symplectic structure on
$\CM_b(\bX,\alpha)$. It is tantalizingly similar to R.~Bielawski's formula,
but in general (rank $>1$) the two structures do not coincide. Let
$\bP\supset\bB$ be a parabolic subgroup. The construction of the Poisson
structure on $\CM_b(\bX,\alpha)$ generalizes {\em verbatim} to the space
of based maps $\CM=\CM_b(\bG/\bP,\beta)$. In most cases the corresponding map
$T^*\CM\to T\CM$ is not an isomorphism, i.e. $\CM$ splits into nontrivial
symplectic leaves. These leaves are explicitly described.
\end{abstract}
\maketitle

\section{Introduction}

\subsection{} Let $\bG$ be a semisimple complex Lie group with the Cartan
datum $(I,\cdot)$ and the root datum $(Y,X,\ldots)$.
Let $\bH\subset\bB=\bB_+,\bB_-\subset\bG$ be a Cartan subgroup and a pair of
opposite Borel subgroups respectively.
Let $\bX=\bG/\bB$ be the flag manifold of $\bG$. Let $C=\PP^1\ni\infty$
be the projective line.
Let $\alpha=\sum_{i\in I}a_i\alpha_i\in{\mathbb N}[I]\subset
H_2(\bX,{\mathbb Z})$.

The moduli space of {\em $\bG$-monopoles of topological charge $\alpha$}
(see e.g. ~\cite{j}) is naturally identified with the space $\CM_b(\bX,\alpha)$
of based maps from $(C,\infty)$ to $(\bX,\bB_+)$ of degree $\alpha$.
The moduli space of
$\bG$-monopoles carries a natural hyperk\"ahler structure, and hence a
holomorphic symplectic structure. We propose a simple explicit formula for
the symplectic structure on $\CM_b(\bX,\alpha)$. It generalizes the
well known formula for $\bG=SL_2$ ~\cite{ah}.

\subsection{} Recall that for $\bG=SL_2$ we have $(\bX,\bB_+)=(\PP^1,\infty)$.
Recall the natural local coordinates on $\CM_b(\PP^1,a)$ (see ~\cite{ah}).
We fix a coordinate $z$ on $C$ such that $z(\infty)=\infty$. Then a based
map $\phi:\ (C,\infty)\to(\PP^1,\infty)$ of degree $a$ is a rational function
$\dfrac{p(z)}{q(z)}$ where $p(z)$ is a degree $a$ polynomial with the leading
coefficient 1, and $q(z)$ is a degree $<a$ polynomial. Let $U$ be the open
subset of based maps such that the roots $x^1,\ldots,x^a$ of $p(z)$ are
multiplicity free. Let $y^k$ be the value of $q(z)$ at $x^k$.
Then $x^1,\ldots,x^a,y^1,\ldots,y^a$ form an \'etale coordinate system on $U$.
The symplectic form $\Omega$ on $\CM_b(\PP^1,a)$ equals
$\sum_{k=1}^a\dfrac{dy^k\wedge dx^k}{y^k}$. In other words, the Poisson brackets
of these coordinates are as follows: $\{x^k,x^m\}=0=\{y^k,y^m\};\
\{x^k,y^m\}=\delta_{km}y^m$.

For an arbitrary $\bG$ we consider the invariant bilinear form $(,)$
on the weight
lattice $X$ such that the square length of a {\em short} simple root
$(\check\alpha_i,\check\alpha_i)=2$. We set
$\check{d}_i:=\frac{(\check\alpha_i,\check\alpha_i)}{2}$.
For $i\in I$ let $\bX_i\subset\bX$ be the corresponding
codimension 1 $\bB_-$-orbit (Schubert cell), and let
$\overline\bX_i\supset\bX_i$
be its closure (Schubert variety). For $\phi\in\CM_b(\bX,\alpha)$ we define
$x_i^1,\ldots,x_i^{a_i}\in\BA^1$ as the points of intersection of $\phi(\PP^1)$
with $\overline\bX_i$. This way we obtain the projection
$\pi^\alpha:\ \CM_b(\bX,\alpha)\to\BA^\alpha$ (the configuration space of
$I$-colored divizors of degree $\alpha$ on $\BA^1$).
Let $U\subset\CM_b(\bX,\alpha)$ be the open subset of based maps such that
$\phi(\PP^1)\cap\overline\bX_i\subset\bX_i$ for any $i\in I$, and
$x_i^k\ne x_j^l$ for any $i,j\in I,\ 1\leq k\leq a_i,\ 1\leq l\leq a_j$.
Locally in $\bX$ the cell $\bX_i$ is the zero divizor of a function $\varphi_i$
(globally, $\varphi_i$ is a section of the line bundle $L_{\check\omega_i}$
corresponding to the fundamental weight $\check\omega_i\in X$). The rational function
$\varphi_i\circ\phi$ on $C$ is of the form $\dfrac{p_i(z)}{q_i(z)}$ where
$p_i(z)$ is a degree $a_i$ polynomial with the leading
coefficient 1, and $q_i(z)$ is a degree $<a_i$ polynomial.
Let $y^k_i$ be the value of $q_i(z)$ at $x^k_i$.
Then $x^k_i,y^k_i,\ i\in I,\ 1\leq k\leq a_i$,
form an \'etale coordinate system on $U$.
The Poisson brackets of these coordinates are as follows:
$$\{x_i^k,x_j^l\}=0=\{y_i^k,y_i^l\};\
\{x_i^k,y_j^l\}=\check{d}_i\delta_{ij}\delta_{kl}y_j^l;\
\{y_i^k,y_j^l\}=(\check\alpha_i,\check\alpha_j)\frac{y_i^ky_j^l}{x_i^k-x_j^l}\
\operatorname{for}\ i\ne j.$$

\subsection{} It follows that the symmetric functions of the $x$-coordinates
(well defined on the whole $\CM_b(\bX,\alpha)$)
are in involution. In other words, the projection
$\pi^\alpha:\ \CM_b(\bX,\alpha)\to\BA^\alpha$ is an integrable system on
$\CM_b(\bX,\alpha)$. The fibers of $\pi^\alpha:\ U\to\BA^\alpha$ are
Lagrangian submanifolds of $U$. It is known that all the fibers of
$\pi^\alpha:\ \CM_b(\bX,\alpha)\to\BA^\alpha$ are equidimensional of the same
dimension $|\alpha|$ (see ~\cite{fm}), hence $\pi^\alpha$ is flat, hence
all the fibers are Lagrangian.

\subsection{} Let $\bP\supset\bB$ be a parabolic subgroup. The construction
of the Poisson structure on $\CM_b(\bX,\alpha)$ generalizes {\em verbatim}
to the space of based maps $\CM=\CM_b(\bG/\bP,\beta)$. In most cases the
corresponding map $P:\ T^*\CM\to T\CM$ is not an isomorphism, i.e. $\CM$ splits
into nontrivial symplectic leaves. For certain degrees $\alpha\in\BN[I]$
we have the natural embedding $\Pi:\ \CM_b(\bX,\alpha)\hookrightarrow\CM$,
and the image is a symplectic leaf of $P$. Moreover, all the symplectic leaves
are of the form $g\Pi(\CM_b(\bX,\alpha))$ for certain $\alpha\in\BN[I],\
g\in\bP$, see the Theorem ~\ref{tautology}.

\subsection{} The above Poisson structure is a baby (rational) version of
the Poisson structure on the moduli space of $\bB$-bundles over an elliptic
curve ~\cite{fo}. We learnt of its definition (as a differential in the
hypercohomology spectral sequence, see \S\ref{second}) from B.Feigin.
Thus, our modest contribution reduces just to a proof of Jacobi identity.
Note that the Poisson structure of ~\cite{fo} arises as a quasiclassical
limit of {\em elliptic algebras}. On the other hand, $\CM_b(\bX,\alpha)$
is an open subset in the moduli space $\hat\CM_b(\bX,\alpha)$ of $\bB$-bundles
on $C$ trivialized at $\infty$, such that the induced $\bH$-bundle has degree
$\alpha$. One can see easily that $\hat\CM_b(\bX,\alpha)$ is isomorphic to
an affine space $\BA^{2|\alpha|}$, and the symplectic structure on
$\CM_b(\bX,\alpha)$ extends to the Poisson structure on $\hat\CM_b(\bX,\alpha)$.
The latter one can be quantized along the lines of ~\cite{fo}.


\subsection{Notations}

For a subset $J\subset I$ we denote by $\bP_J\supset\bB$ the corresponding
parabolic subgroup. Thus, $\bP_\emptyset=\bB$.
Denote by $\bX_J=\bG/{\bP_J}$ the corresponding
parabolic flag variety; thus, $\bX_\emptyset=\bX$. We denote by
$\varpi:\ \bX\to\bX_J$ the natural projection.
We denote by $x\in \bX_J$ the marked point $\varpi(\bB_+)$.

Let $\CM=\CM_b(\bX_J,\al)$ denote the space of based
algebraic maps $\phi:(C,\infty)\to(\bX_J,x)$ of degree
$\al\in H_2(\bX_J,\ZZ)$.

Let $\fg$ denote the Lie algebra of $\bG$.
Let $\fg_{\bX_J}$ denote the trivial vector bundle with
the fiber $\fg$ over ${\bX_J}$ and let $\fp_{\bX_J}\subset\fg_{\bX_J}$
(resp. $\fr_{\bX_J}\subset\fp_{\bX_J}\subset\fg_{\bX_J}$) be its subbundle
with the fiber over a point $P$ equal to the corresponding
Lie subalgebra $\fp\subset\fg$ (resp. its nilpotent
radical $\fr\subset\fp\subset\fg$). In case $J=\emptyset$ we will also
denote $\fr_{\bX_\emptyset}$ by $\fn_\bX$, and $\fp_{\bX_\emptyset}$ by
$\fb_\bX$. Note that the quotient bundle $\fh_\bX:=\fb_\bX/\fn_\bX$ is trivial
(abstract Cartan algebra).

Recall that the tangent bundle $T{\bX_J}$ of ${\bX_J}$ (resp. cotangent
bundle $T^*{\bX_J}$) is canonicaly isomorphic to the bundle
$\fg_{\bX_J}/\fp_{\bX_J}$ (resp. $\fr_{\bX_J}$).

\subsection{Acknowledgments.}

This paper has been written during the stay
of the second author at the Max-Planck-Institut
f\"ur Mathematik. He would like
to express his sincere gratitude to the Institut
for the hospitality and the excellent work conditions.
It is clear from the above discussion that the present note owes its
existense to the generous explanations of B.Feigin. We are deeply grateful
to L.~Rybnikov who has found and corrected a mistake in the calculation of
the symplectic form in the non simply laced case.

\section{The Poisson structure}\label{second}

\subsection{}\label{def}

The fibers of the tangent and cotangent
bundles of the space $\CM$ at the point $\phi$ are computed as follows:
$$
\arraycolsep=1pt
\begin{array}{lll}
T_\phi\CM&=H^0(C,(\phi^*T{\bX_J})\om)&=
H^0(C,(\phi^*\fg_{\bX_J}/\fp_{\bX_J})\om),\\
T^*_\phi\CM&=H^1(C,(\phi^*T^*{\bX_J})\om)&=H^1(C,(\phi^*\fr_{\bX_J})\om).
\end{array}
$$
The second identification follows from the first by the Serre
duality.

We have a tautological complex of vector bundles on $\bX_J$:
\begin{equation}\label{nconv}
\fr_{\bX_J} \to \fg_{\bX_J} \to \fg_{\bX_J}/\fp_{\bX_J}
\end{equation}
The pull-back via $\phi$ of this complex twisted by $\CO_C(-1)$
gives the following complex of vector bundles on $C$
\begin{equation}\label{main}
(\phi^*\fr_{\bX_J})\om \to (\phi^*\fg_{\bX_J})\om \to
(\phi^*\fg_{\bX_J}/\fp_{\bX_J})\om
\end{equation}

Consider the hypercohomology spectral sequence of the
complex \refeq{main}. Since $\fg_{\bX_J}$ is the trivial vector
bundle we have $H^\bullet(C,(\phi^*\fg_{\bX_J})\om)=0$, hence the
second differential of the spectral sequence induces a map
$$
d_2:H^1(C,(\phi^*\fr_{\bX_J})\om)\to H^0(C,(\phi^*\fg_{\bX_J}/\fp_{\bX_J})\om)
$$
that is a map $P^{\bX_J}_\phi:T^*_\phi\CM\to T_\phi\CM$.
This construction easily globalizes to give a morphism
$P^{\bX_J}:\ T^*\CM\to T\CM$.

\begin{thm}\label{poisson}
$P$ defines a Poisson structure on $\CM$.
\end{thm}

Here we will reduce the Theorem to the case $\bX_J=\bX$.
This case will be treated in the next section.

\subsection{}\label{reduction}

Let $\varpi_*:H_2(\bX,\ZZ)\to H_2({\bX_J},\ZZ)$ be
the push-forward map. The map $\varpi$ induces a map
$\Pi:\CM_b(\bX,\al)\to\CM_b({\bX_J},\varpi_*\al)$.

\begin{pr}\label{F}
The map $\Pi$ respects $P$, that is the following square is commutative
$$
\begin{CD}
T^*\CM_b(\bX,\al)           @>{P^\bX}>> T\CM_b(\bX,\al)       \\
@A{\Pi^*}AA                         @V{\Pi_*}VV         \\
\Pi^*T^*\CM_b({\bX_J},\varpi_*\al) @>{P^{\bX_J}}>>
\Pi^*T\CM_b({\bX_J},\varpi_*\al)
\end{CD}
$$
\end{pr}
\begin{proof}
We have the following commutative square on $\bX$
$$
\begin{CD}
\fn_\bX      @>>> \fg_\bX      @>>> \fg_\bX/\fb_\bX     \\
@AAA		@|		@|		\\
\varpi^*\fr_{\bX_J}      @>>> \fg_\bX      @>>> \fg_\bX/\fb_\bX     \\
@|		@|		@VVV		\\
\varpi^*\fr_{\bX_J} @>>> \varpi^*\fg_{\bX_J} @>>>
\varpi^*(\fg_{\bX_J}/\fp_{\bX_J})
\end{CD}
$$
Consider its pull-back via $\phi\in\CM_b(\bX,\al)$ twisted by $\CO_C(-1)$.
Let $d_2$ denote the second differential of the hypercohomology spectral
sequence of the middle row. Then we evidently have
$$
P^\bX\cdot\Pi^*=d_2,\qquad\Pi_*\cdot d_2=P^{\bX_J}
$$
and the Proposition follows.
%
%
\end{proof}

Now, assume that we have proved that $P^\bX$ defines
a Poisson structure. For any $\beta\in H_2({\bX_J},\ZZ)$
we can choose $\al\in H_2(\bX,\ZZ)$ such that $\varpi_*\al=\beta$
and the map $\Pi$ is open. Then the algebra of
functions on $\CM_b({\bX_J},\beta)$ is embedded into the algebra
of functions on $\CM_b(\bX,\al)$ and the Proposition \ref{F}
shows that the former bracket is induced by the latter one. Hence it
is also a Poisson bracket.


\section{The case of $\bX$}

In this section we will denote $\CM_b(\bX,\al)$ simply by $\CM$.

\subsection{}

Since $\fh$ is a trivial vector bundle on $\bX$ the exact
sequences
$$
0\to\fn_\bX\to\fb_\bX\to\fh_\bX\to0,\qquad
0\to\fh_\bX\to\fg_\bX/\fn_\bX\to\fg_\bX/\fb_\bX\to0
$$
induce the isomorphisms
$$
\arraycolsep=1pt
\begin{array}{lll}
T_\phi\CM & =H^0(C,(\phi^*\fg_\bX/\fb_\bX)\om) & =
H^0(C,(\phi^*\fg_\bX/\fn_\bX)\om),\\
T^*_\phi\CM & = H^1(C,(\phi^*\fn_\bX)\om) & = H^1(C,(\phi^*\fb_\bX)\om)
\end{array}
$$
Applying the construction of \ref{def} to the
following tautological complex of vector bundles on $\bX$
\begin{equation}\label{conv}
\fb_\bX\to\fg_\bX\oplus\fh_\bX\to\fg_\bX/\fn_\bX
\end{equation}
and taking into account the above isomorphisms we
get a map $\tilde P^\bX_\phi:T^*_\phi\CM\to T_\phi\CM$.

\begin{lem} We have $\tilde P^\bX_\phi=P^\bX_\phi$.
\end{lem}
\begin{proof}
The same reasons as in the proof of the Proposition \ref{F} work
if we consider the following commutative diagram
$$
\begin{CD}
\fn_\bX @>>> \fg_\bX            @>>> \fg_\bX/\fb_\bX    \\
@|           @| 			 @AAA		\\
\fn_\bX @>>> \fg_\bX            @>>> \fg_\bX/\fn_\bX    \\
@VVV       @VVV			 @|		\\
\fb_\bX @>>> \fg_\bX\oplus\fh_\bX @>>> \fg_\bX/\fn_\bX
\end{CD}
$$
\end{proof}

It will be convenient for us to use the complex \refeq{conv}
for the definition of the map $P^\bX_\phi$ instead of \refeq{nconv}.

\subsection{}

Here we will describe the Pl\"ucker embedding of the space $\CM$.

Let $X\supset\mho\cong I$ be the set of fundamental weights:
$\langle\alpha_i,\check\omega_j\rangle=\delta_{ij}$.
Recall that we consider the invariant bilinear form $(,)$ on the weight
lattice $X$ such that the square length of a {\em short} simple root
$(\check\alpha_i,\check\alpha_i)=2$. We set
$\check{d}_i:=\frac{(\check\alpha_i,\check\alpha_i)}{2}$.
For a dominant weight ${\check\lambda}\in X$
we denote by $V_{\check\lambda}$ the irreducible $G$-module with highest weight
${\check\lambda}$.

Recall that $\bX$ is canonically embedded into the
product of projective spaces
$$
\bX\subset\prod_{\o\in\Om}\PP(V_\o)
$$
This induces the embedding
$$
\CM\subset\prod_{\o\in\Om}\CM_b(\PP(V_\o),\langle\alpha,\check\omega\rangle)
$$
Note that the marked point of the space $\PP(V_\o)$
is just the highest weight vector $v_\o$ with respect to
the Borel subgroup $\bB$.

A degree $d$ based map $\phi_\o:(C,\infty)\to(\PP(V_\o),v_\o)$
can be represented by a $V_\o$-valued degree $d$ polynomial in $z$,
taking the value $v_\o$ at infinity. Let us denote the affine
space of such polynomials by $R_d(V_\o)$.

The {\em Pl\"ucker embedding} of the space $\CM$ is the
embedding into the product of affine spaces
$$
\CM\subset\prod_{\o\in\Om}R_{\langle\alpha,\check\omega\rangle}(V_\o).
$$


A map $\phi\in\CM$ will be represented by a collection of
polynomials
$(\phi_\o\in R_{\langle\alpha,\check\omega\rangle}(V_\o))_{\o\in\Om}$.

\subsection{The coordinates}

The dual representation $V^*_\o$ decomposes into
the sum of weight subspaces
$$
V^*_\o=\oplus_{\la\in X}{V^*}_\o^\la.
$$

We choose a weight base $(f_\o^\la)$ of $V_\o^*$,
such that $f_\o^{-\o}(v_\o)=1$.
Suppose $\langle\alpha_i,\check\omega\rangle=1$. Then $\dim{V^*}_{\check\omega}^{-{\check\omega}}=
\dim{V^*}_{\check\omega}^{-{\check\omega}+\check\alpha_i}=1$, and
$\dim{V^*}_{\check\omega}^{-{\check\omega}+\check\alpha_i+\check\alpha_j}=0$ if $(\check\alpha_i,\check\alpha_j)=0$, and
$\dim{V^*}_{\check\omega}^{-{\check\omega}+\check\alpha_i+\check\alpha_j}=1$ if $(\check\alpha_i,\check\alpha_j)\ne0$. Hence, in the latter case,
the vectors $f_\o^{-\o}$, $f_\o^{\check\alpha_i-\o}$ and
$f_\o^{\check\alpha_i+\check\alpha_j-\o}$ are defined uniquely up to multiplication by a constant.
Let $E_i,F_i,H_i$ be the standard generators of $\fg$. Then we will take
$f^{\check\alpha_i-{\check\omega}}_{\check\omega}:=E_if^{-{\check\omega}}_{\check\omega},\ f^{\check\alpha_i+\check\alpha_j-{\check\omega}}_{\check\omega}:=
E_jE_if^{-{\check\omega}}_{\check\omega}$.

We consider the polynomials $\phi_\o^\la:=f_\o^\la(\phi_\o)$:
the $\la$ weight components of $\phi_\o$. In particular,
$\phi_\o^{-\o}$ is the degree $\langle\alpha,{\check\omega}\rangle$ unitary polynomial
and $\phi_\o^{\check\alpha_i-\o}$ is the degree $<\langle\alpha,{\check\omega}\rangle$ polynomial.

Let $x_\o^1,\dots,x_\o^{\langle\alpha,{\check\omega}\rangle}$
be the roots of $\phi_\o^{-\o}$
and $y_\o^1,\dots,y_\o^{\langle\alpha,{\check\omega}\rangle}$
be the values of $\phi_\o^{\check\alpha_i-\o}$
at the points $x_\o^1,\dots,x_\o^{\langle\alpha,{\check\omega}\rangle}$ respectively.
Consider the open subset $U\subset\CM$ formed by all the
maps $\phi$ such that all $x_\o^k$ are distinct
and all $y_\o^k$ are non-zero.
On this open set we have
$$
\phi_\o^{-\o}(z)=\prod_{k=1}^{\langle\alpha,{\check\omega}\rangle}(z-x_\o^k),\qquad
\phi_\o^{\check\alpha_i-\o}(z)=\sum_{k=1}^{\langle\alpha,{\check\omega}\rangle}
\frac{y_\o^k\phi_\o^{-\o}(z)}{(\phi_\o^{-\o})'(x_\o^k)(z-x_\o^k)}.
$$
The collection of $2|\al|$ functions
\begin{equation}\label{xy}
(x_\o^k,y_\o^k),\qquad(\o\in\Om,\ 1\le k\le\langle\alpha,{\check\omega}\rangle)
\end{equation}
is an \'etale coordinate system in $U$.
One can either check this straightforwardly,
or just note that the matrix of $P^\bX$ in these coordinates has a maximal
rank, see the Remark ~\ref{matrix} below.

So let us compute the map $P^\bX$ in these coordinates.

\subsection{}

We have an isomorphism $V_\o^*\cong V_{\o^*}$ for a certain involution
$\o\mapsto\o^*$ of $\Om$. Note that the induced involution of the weight
lattice $X$ preserves the invariant bilinear form $(,)$.
The action of $\fg$ on $V_{\o^*}$ induces an embedding
$$
\fg_\bX/\fn_\bX\subset\opl_{\o\in\Om}V_{\o^*}\otimes L_\o
$$
of vector bundles over $\bX$ and the dual surjection
$$
\opl_{\o\in\Om}V^*_{\o^*}\otimes L^*_\o\to\fb_\bX,
$$
where $L_\o$ stands for the line bundle, corresponding
to the weight $\o$.
Hence we have the following complex
\begin{equation}\label{newc}
\opl_{\o\in\Om}V^*_{\o^*}\otimes L^*_\o\to\fg_\bX\oplus\fh_\bX
\to\opl_{\o\in\Om}V_{\o^*}\otimes L_\o
\end{equation}

\begin{rem}
The differentials in the above complex in the fiber over a
point $\bB'\in\bX$
are computed as follows:
$$
\varphi\in V^*_\o\mapsto
\left(\sum\varphi(\xi^kv')\xi_k\right)\oplus
\left(\sum\o(h^i)\varphi(v')h_i\right)\in\fg\oplus\fh,
$$
$$
\xi\oplus h\in\fg\oplus\fh\mapsto\xi v'-\o(h)v'\in V_\o.
$$
Here $v'$ is a highest weight vector of $V_\o$ with respect
to $\bB'$; $(\xi_k)$, $(\xi^k)$ are dual (with respect to the
standard scalar product) bases of $\fg$; and $(h_i)$, $(h^i)$
are dual bases of $\fh$.
\end{rem}

\nc{\CF}{{\mathcal F}}
\nc{\CG}{{\mathcal G}}
\nc{\tr}{\operatorname{tr}\nolimits}

\subsection{}

In order to compute the brackets of the coordinates \refeq{xy}
at a point $\phi\in\CM$ we need to take the pull-back of the
complex \refeq{newc} via $\phi$, twist it by $\CO_C(-1)$
and compute the second differential of the hypercohomology spectral
sequence. The following Lemma describes this differential
in general situation.

\begin{lem}\label{D}
Consider a complex
$\begin{CD}
K^\bullet=(\CF @>{\ f\ }>> A\otimes\CO_C @>{\ g\ }>> \CG)
\end{CD}$
on $C$, where $A$ is a vector space and
$$
\begin{array}{rclcl}
f & \in & Hom(\CF,A\otimes\CO_C) & = & A\otimes H^0(C,\CF^*),\\
g & \in & Hom(A\otimes\CO_C,\CG) & = & A^*\otimes H^0(C,\CG).
\end{array}
$$
Consider
$$
D=\tr(f\otimes g)\in H^0(C,\CF^*)\otimes H^0(C,\CG)=
H^0(C\times C,\CF^*\boxtimes\CG),
$$
where $\tr:A\otimes A^*\to\CCC$ is the trace homomorphism.
Then

$1)$ The restriction of $D$ to the diagonal $\Delta\subset C\times C$
vanishes, hence $D=\tilde D\Delta$ for some
\begin{multline*}
\tilde D\in H^0(C\times C,(\CF^*\boxtimes\CG)(-\Delta))=
H^0(C,\CF^*(-1))\otimes H^0(C,\CG(-1))=\\=
H^1(C,\CF(-1))^*\otimes H^0(C,\CG(-1)).
\end{multline*}

$2)$ The second differential $d_2:H^1(C,\CF(-1))\to H^0(C,\CG(-1))$ of
the hypercohomology spectral sequence of $K^\bullet\otimes\CO_C(-1)$
is induced by the section $\tilde D$.
\end{lem}
\begin{proof}
The first statement is evident. To prove the second statement
consider the following commutative diagram on $C\times C$
$$
\begin{CD}
\CF(-1)\boxtimes\CO_C			\!@>{f(-1)\boxtimes1}>>
A\otimes\CO_{C\times C}(-1,0)           @>{(1\boxtimes g)|_{\Delta}}>>
(\CO_C(-1)\boxtimes\CG)|_{\Delta}       \\
@V{\tilde D}VV		@V{1\boxtimes g}VV	@|	\\
\CO_C(-2)\boxtimes\CG(-1)               \!@>\Delta>>
\CO_C(-1)\boxtimes\CG                   @>{|_\Delta}>>
(\CO_C(-1)\boxtimes\CG)|_{\Delta}
\end{CD}
$$
Both rows are complexes with acyclic middle term, hence the
second differentials of the hypercohomology spectral sequences
commute with the maps induced on cohomology by the vertical arrows:
$$
\begin{CD}
H^1(C\times C,\CF(-1)\boxtimes\CO_C)			@>d_2>>
H^0(C\times C,(\CO_C(-1)\boxtimes\CG)|_{\Delta})        \\
@V{\tilde D}VV				@|		\\
H^1(C\times C,\CO_C(-2)\boxtimes\CG(-1))                @>d_2>>
H^0(C\times C,(\CO_C(-1)\boxtimes\CG)|_{\Delta})
\end{CD}
$$
Now it remains to note that
$$
\begin{array}{lcl}
H^1(C\times C,\CF(-1)\boxtimes\CO_C) & = & H^1(C,\CF(-1)),		\\
H^1(C\times C,\CO_C(-2)\boxtimes\CG(-1)) & = & H^0(C,\CG(-1)),          \\
H^0(C\times C,(\CO_C(-1)\boxtimes\CG)|_{\Delta}) & = & H^0(C,\CG(-1)),
\end{array}
$$
and that the map $H^0(C,\CG(-1))\to H^0(C,\CG(-1))$ induced
by the map $d_2$ in the second row of the above diagram is identity.
\end{proof}

\subsection{}

Consider the pullback of \refeq{newc} via $\phi\in\CM$, and
twist it by $\CO_C(-1)$. We want to apply Lemma \ref{D} to compute
the $(\o_i,{\check\omega}_j)$-component of the second differential of the
hypercohomology spectral sequence.

In notations of the Lemma we have
\begin{multline*}
D_{{\check\omega}_i,{\check\omega}_j}(z,w)=
\sum\xi^k\phi_{{\check\omega}_i}(z)\otimes\xi_k\phi_{{\check\omega}_j}(w)-
\sum{{\check\omega}_i^*}(h^i)\phi_{{\check\omega}_i}(z)\otimes
{\check\omega}_j^*(h_i)\phi_{{\check\omega}_j}(w)=\\=
\sum\xi^k\phi_{{\check\omega}_i}(z)\otimes\xi_k\phi_{{\check\omega}_j}(w)-
({{\check\omega}_i},{\check\omega}_j)\phi_{{\check\omega}_i}(z)\otimes\phi_{{\check\omega}_j}(w)\in
V_{{\check\omega}_i^*}\otimes V_{{\check\omega}_j^*}(z,w).
\end{multline*}

\begin{lem}\label{lie}
The operator $\sum\xi^k\otimes\xi_k-(\o_i,\o_j)$ acts as
a scalar multiplication on every irreducible summand
$V_\la\subset V_{{\check\omega}_i}\otimes V_{{\check\omega}_j}$.

On $V_{{{\check\omega}_i}+{\check\omega}_j}$ it acts as a multiplication by $0$.

If ${{\check\omega}_i}={\check\omega}_j$, then on $V_{2{{\check\omega}_i}-\check\alpha_i}\subset V_{{\check\omega}_i}
\otimes V_{{\check\omega}_i}$ it acts as a multiplication by
$-(\check\alpha_i,\check\alpha_i)=-2\check{d}_i$.

If $i\ne j$, $(\check\alpha_i,\check\alpha_j)\ne0$ then on
$V_{{{\check\omega}_i}+{\check\omega}_j-\check\alpha_i-\check\alpha_j}\subset V_{{\check\omega}_i}\otimes V_{{\check\omega}_j}$ it acts
as a multiplication by $((\check\alpha_i,\check\alpha_j)-
\frac{(\check\alpha_i,\check\alpha_i)}{2}-
\frac{(\check\alpha_j,\check\alpha_j)}{2})$.
\end{lem}
\begin{proof}
It is easy to check that $\sum\xi^k\otimes\xi_k$ commutes
with the natural action of $\fg$ on $V_{{\check\omega}_i}\otimes V_{{\check\omega}_j}$.
The first part of the Lemma follows. The rest of the Lemma
can be checked by the straightforward computation of the action
of $\sum\xi^k\otimes\xi_k$ on the highest vectors of the
corresponding subrepresentations.
\end{proof}

\subsection{}

If we want to compute the brackets of the coordinates \refeq{xy}
we are interested in the components of $D_{{{\check\omega}_i^*},{\check\omega}_j^*}(z,w)$
in the weights
\begin{equation}\label{weights}
{{\check\omega}_i}+{\check\omega}_j,\quad{{\check\omega}_i}+{\check\omega}_j-\check\alpha_i,\quad{{\check\omega}_i}
+{\check\omega}_j-\check\alpha_j,\quad{{\check\omega}_i}+{\check\omega}_j-\check\alpha_i-\check\alpha_j.
\end{equation}
The following Lemma describes the corresponding weight components
of the tensor product $V_{{\check\omega}_i}\otimes V_{{\check\omega}_j}$.
\begin{lem}\label{lie1}
The embedding $V_{{{\check\omega}_i}+{\check\omega}_j}\subset V_{{\check\omega}_i}\otimes V_{{\check\omega}_j}$ induces
an isomorphism in the weights \refeq{weights} with the following
two exceptions:

\noindent
$(1)$ $(V_{{\check\omega}_i}\otimes V_{{\check\omega}_i})^{2{{\check\omega}_i}-\check\alpha_i}=
V_{2{{\check\omega}_i}}^{2{{\check\omega}_i}-\check\alpha_i}\oplus V_{2{{\check\omega}_i}-\check\alpha_i}^{2{{\check\omega}_i}-\check\alpha_i}$;
the $G$-projection
to the second summand is given by the formula
$$
a(v_{\o_i}\otimes F_iv_{\o_i})+b(F_iv_{\o_i}\otimes v_{\o_i})
\mapsto\frac{a-b}2(v_{\o_i}\otimes F_iv_{\o_i}-F_iv_{\o_i}\otimes v_{\o_i}).
$$

\noindent
$(2)$
$(V_{{\check\omega}_i}\otimes V_{{\check\omega}_j})^{{{\check\omega}_i}+{\check\omega}_j-\check\alpha_i-\check\alpha_j}=
V_{{{\check\omega}_i}+{\check\omega}_j}^{{{\check\omega}_i}+{\check\omega}_j-\check\alpha_i-\check\alpha_j}\oplus
V_{{{\check\omega}_i}+{\check\omega}_j-\check\alpha_i-\check\alpha_j}^{{{\check\omega}_i}+{\check\omega}_j-\check\alpha_i-\check\alpha_j}$
if $i\ne j$ and $(\check\alpha_i,\check\alpha_j)\ne0$;
the $G$-projection to the second summand is given by the formula
\begin{multline*}
a(v_{\o_i}\otimes F_iF_jv_{\o_j})+b(F_iv_{\o_i}\otimes F_jv_{\o_j})+
c(F_jF_iv_{\o_i}\otimes v_{\o_j})\mapsto\\\mapsto
\frac{b-a-c}{1-\langle\alpha_i,\check\alpha_j\rangle^{-1}-
\langle\alpha_j,\check\alpha_i\rangle^{-1}}
(\langle\alpha_i,\check\alpha_j\rangle^{-1}v_{{\check\omega}_i}\otimes F_iF_jv_{{\check\omega}_j}+
F_iv_{{\check\omega}_i}\otimes F_jv_{\o_j}+
\langle\alpha_j,\check\alpha_i\rangle^{-1}F_jF_iv_{{\check\omega}_i}\otimes v_{{\check\omega}_j}).
\end{multline*}
\end{lem}
\begin{proof}
Straightforward.
\end{proof}

\subsection{}\label{td}

Hence (see Lemma \ref{lie}, Lemma \ref{lie1}) when $\la^*$ is one
of the weights \refeq{weights} the $\la^*$-component
$\tilde D_{{{\check\omega}_i},{\check\omega}_j}^\la(z,w)$ of the polynomial
$\tilde D_{{{\check\omega}_i},{\check\omega}_j}(z,w)=\dfrac{D_{{{\check\omega}_i},{\check\omega}_j}(z,w)}{z-w}$
is zero with the following two exceptions

\begin{equation}\label{ooi}
\tilde D_{{{\check\omega}_i},{{\check\omega}_i}}^{2{{\check\omega}_i}-\check\alpha_i}=
\check{d}_i\frac{\phi_{{\check\omega}_i}^{{\check\omega}_i}(z)\phi_{{\check\omega}_i}^{{{\check\omega}_i}-\check\alpha_i}(w)-
\phi_{{\check\omega}_i}^{{{\check\omega}_i}-\check\alpha_i}(z)\phi_{{\check\omega}_i}^{{\check\omega}_i}(w)}{z-w}
(F_{i^*}v_{\o_i^*}\otimes v_{\o_i^*}-v_{\o_i^*}\otimes F_{i^*}v_{\o_i^*})
\end{equation}
\begin{multline}\label{ooii}
\tilde D_{{{\check\omega}_i},{\check\omega}_j}^{{{\check\omega}_i}+{\check\omega}_j-\check\alpha_i-\check\alpha_j}\!=\!
\frac{\phi_{{\check\omega}_i}^{{{\check\omega}_i}-\check\alpha_i}(z)\phi_{{\check\omega}_j}^{{\check\omega}_j-\check\alpha_j}(w)-
\phi_{{\check\omega}_i}^{{\check\omega}_i}(z)\phi_{{\check\omega}_j}^{{\check\omega}_j-\check\alpha_i-\check\alpha_j}(w)-
\phi_{{\check\omega}_i}^{{{\check\omega}_i}-\check\alpha_i-\check\alpha_j}(z)\phi_{{\check\omega}_j}^{{\check\omega}_j}(w)}{z-w}\times\\
\times(\check\alpha_i,\check\alpha_j)
(\langle\alpha_i,\check\alpha_j\rangle^{-1}v_{{\check\omega}_i^*}\otimes F_{i^*}F_{j^*}v_{{\check\omega}_j^*}+
F_{i^*}v_{{\check\omega}_i^*}\otimes F_{j^*}v_{\o_j^*}+
\langle\alpha_j,\check\alpha_i\rangle^{-1}F_{j^*}F_{i^*}v_{{\check\omega}_i^*}\otimes v_{{\check\omega}_j^*})
\end{multline}

Note that the scalar multiplicators of~Lemma~\ref{lie} have almost canceled
the denominators of~Lemma~\ref{lie1}.

\subsection{}

Now we can compute the brackets.

\begin{pr}\label{formulas}
We have
\begin{equation}\label{forms}
\arraycolsep=0pt
\def\arraystretch{2}
\begin{array}{crcrccrl}
\{&x_{{\check\omega}_i}^k&,&x_{{\check\omega}_j}^l&\} &\ =\ & 0;\\
\{&x_{{\check\omega}_i}^k&,&y_{{\check\omega}_j}^l&\} &\ =\ &
\check{d}_i\delta_{kl}\delta_{ij}y_{{\check\omega}_j}^l;\\
\{&y_{{\check\omega}_i}^k&,&x_{{\check\omega}_j}^l&\} &\ =\ &
-\check{d}_i\delta_{kl}\delta_{ij}y_{{\check\omega}_i}^k;\\
\{&y_{{\check\omega}_i}^k&,&y_{{\check\omega}_j}^l&\} &\ =\ &
(\check\alpha_i,\check\alpha_j)\dfrac{y_{{\check\omega}_i}^ky_{{\check\omega}_j}^l}{x_{{\check\omega}_i}^k-x_{{\check\omega}_j}^l}, &
\qquad\text{if $i\ne j$;}\\
\{&y_{{\check\omega}_i}^k&,&y_{{\check\omega}_i}^l&\} &\ =\ & 0.
\end{array}
\end{equation}
\end{pr}
\begin{proof}
Note that if $p\in V_{\o_i}(z)$ then
$$
dy_{\o_i}^k(p)=\left\langle f_{\o_i}^{\check\alpha_i-\o_i}, p(x_{\o_i}^k)\right\rangle,\qquad
dx_{\o_i}^k(p)=\left\langle f_{\o_i}^{-\o_i},
\frac{p(x_{\o_i}^k)}{(\phi_{\o_i}^{\o_i})'(x_{\o_i}^k)}\right\rangle,
$$
where $\langle\bullet,\bullet\rangle$ stands for the natural pairing.
Note also that
$$
\phi_{{\check\omega}_i}^{{\check\omega}_i}(x_{{\check\omega}_i}^k)=0,\qquad
\phi_{{\check\omega}_i}^{{{\check\omega}_i}-\check\alpha_i}(x_{{\check\omega}_i}^k)=y_{{\check\omega}_i}^k
$$
by definition and
$$
\langle f_{\o_i}^{\check\alpha_i-\o_i},F_iv_{\o_i}\rangle=
\langle E_if_{\o_i}^{-\o_i},F_iv_{\o_i}\rangle=-
\langle f_{\o_i}^{-\o_i},E_iF_iv_{\o_i}\rangle=-1.
$$
Now the Proposition follows from the Lemma \ref{D} and from
the formulas of \ref{td}.
\end{proof}

\begin{rem}\label{matrix}
The matrix of the bivector field $P^\bX$ in the coordinates
$(x_{\o_i}^k,y_{\o_i}^k)$ looks as follows
$$
\nc{\diag}{\operatorname{diag}}
\def\arraystretch{1.5}
\left(
\begin{array}{c|c}
0 & \diag(\check{d}_iy_{\o_i}^k) \\
\hline
-\diag(\check{d}_iy_{\o_i}^k) & \ast
\end{array}
\right)
$$
Since on the open set $U$ this matrix is evidently nondegenerate
it follows that the functions $(x_{\o_i}^k,y_{\o_i}^k)$ indeed
form an \'etale coordinate system.
\end{rem}

\subsection{}

Now we can prove Theorem \ref{poisson}.

{\em Proof of the Theorem \ref{poisson}}
The reduction to the case $J=\emptyset$ has been done in \ref{reduction}.
The latter case is straightforward by the virtue of
Proposition ~\ref{formulas}. \qed

\begin{cor}
\label{vicious cycle}
The map $P^\bX$ provides the space $\CM_b(\bX,\al)$ with
a holomorphic symplectic structure.
\end{cor}
\begin{proof}
Since $P^\bX$ gives a Poisson structure it suffices to
check that $P^\bX$ is nondegenrate at any point.
To this end recall that the hypercohomology spectral
sequence of a complex $K^\bullet$ converges to
$H^\bullet(\bX,K^\bullet)$. Since the only nontrivial
cohomology of the complex \refeq{nconv} is $\fh_\bX$ in degree
zero, the complex
\refeq{main} is quasiisomorphic to $(\phi^*\fh)\om$ in degree zero, hence
the hypercohomology sequence of the complex \refeq{main}
converges to zero, hence $P^\bX_\phi$ is an isomorphism.
\end{proof}

\begin{rem}
One can easily write down the corresponding
symplectic form  in the coordinates \refeq{xy}:
$$\sum_{i,k}\frac{dy_{\o_i}^k\wedge
dx_{\o_i}^k}{\check{d}_iy_{\o_i}^k}+\sum_{i\ne j}\sum_{k,l}\frac{(\check\alpha_i,\check\alpha_j)}{\check{d}_i\check{d}_j}
\frac{dx_{\o_i}^k\wedge dx_{\o_j}^l}{x_{\o_i}^k-x_{\o_j}^l}$$
\end{rem}

\section{Symplectic leaves}

\subsection{}
We fix $\beta\in\BN[I-J]\subset\ZZ[I-J]=H_2(\bX_J,\ZZ)$, and consider
the Poisson structure on $\CM=\CM_b(\bX_J,\beta)$. In this section we will
describe the symplectic leaves of this structure.

Consider $\alpha\in\BN[I]\subset\ZZ[I]=H_2(\bX,\ZZ)$ such that
$\varpi_*\alpha=\beta$ (see ~\ref{reduction}). Note that $\varpi_*$ is nothing
but the natural projection from $\BN[I]$ to $\BN[I-J]$.
Thus $\alpha-\varpi_*\alpha\in\BN[J]$. We will call an element
$\gamma\in\BN[I]$ {\em $J$-antidominant} if $\langle\gamma,\check\alpha_j\rangle\leq0$ for
any $j\in J$. We will call $\alpha\in\BN[I]$ {\em a special lift} of
$\beta\in\BN[I-J]$ if $\varpi_*\alpha=\beta$, and $\alpha$ is $J$-antidominant.

\begin{lem} If $\alpha$ is a special lift of $\beta$, then the natural
projection $\Pi:\ \CM_b(\bX,\alpha)\to\CM_b(\bX_J,\beta)$ (see ~\ref{reduction})
is an immersion.
\end{lem}

\begin{proof} Let $\phi\in\CM_b(\bX,\alpha)$. Then
$T_\phi\CM_b(\bX,\alpha)=H^0(C,(\phi^*\fg_\bX/\fb_\bX)\otimes\CO_C(-1))$,
and
$T_{\Pi\phi}\CM_b(\bX_J,\beta)=
H^0(C,((\Pi\phi)^*\fg_{\bX_J}/\fp_{\bX_J})\otimes\CO_C(-1))$. Hence the
kernel of the natural map $\Pi_*:\ T_\phi\CM_b(\bX,\alpha)\to
T_{\Pi\phi}\CM_b(\bX_J,\beta)$ equals
$H^0(C,\phi^*(\varpi^*\fp_{\bX_J}/\fb_\bX)\otimes\CO_C(-1))$. Now
$\varpi^*\fp_{\bX_J}/\fb_\bX$ has a natural filtration with the successive
quotients of the form $L_\theta$ where $\theta$ is a positive root of the
root subsystem spanned by $J\subset I$. Since $\alpha$ is a special lift,
$\deg\phi^*L_\theta\leq0$. We conclude that
$H^0(C,(\phi^*L_\theta)\otimes\CO_C(-1))=0$, and thus
$H^0(C,\phi^*(\varpi^*\fp_{\bX_J}/\fb_\bX)\otimes\CO_C(-1))=0$. The Lemma is
proved.
\end{proof}

\begin{rem} For fixed $\beta\in\BN[I-J]$ the set of its special lifts is
evidently finite. It is nonempty (see e.g. the proof
of~Theorem~\ref{tautology}).
\end{rem}

\subsection{}
It follows from the Proposition ~\ref{F} that
$\Pi(\CM_b(\bX,\alpha))$ is a symplectic leaf of the Poisson
structure $P$ on $\CM$, if $\alpha$ is a special lift of $\beta$.

The group $\bP_J$ acts naturally on $\CM$; it preserves $P$ since the
complex ~(\ref{main}) is $\bP_J$-equivariant. It follows that
for $g\in\bP_J$ the subvariety $g\Pi(\CM_b(\bX,\alpha))\subset\CM$ is
also a symplectic leaf. Certainly, $g\Pi(\CM_b(\bX,\alpha))=
\Pi(\CM_b(\bX,\alpha))$ whenever $g\in\bB$.

\begin{thm}
\label{tautology}
Any symplectic leaf of $P$ is of the form
$g\Pi(\CM_b(\bX,\alpha))$ where $\alpha$ is a special lift of $\beta$,
and $g\in\bP_J$.
\end{thm}

\begin{proof}
We only need to check that for any $\psi\in\CM$ there exists a special lift
$\alpha$, a point $\phi\in\CM_b(\bX,\alpha)$, and $g\in\bP_J$ such that
$\psi=g\Pi\phi$. In other words, it suffices to find a special lift
$\alpha$, and a point $\phi\in\CM(\bX,\alpha)$ (unbased maps!) such that
$\phi(\infty)\in\bP_Jx$ (the smallest $\bP_J$-orbit in $\bX$), and
$\psi=\Pi\phi$.

Recall that given a reductive group $G$ with a Cartan subgroup $H$ and
a set of simple roots $\Delta\subset X(H)$,
the isomorphism classes of $G$-torsors over $C$ are numbered
by the set $X^+_*(H)$ of the dominant coweights of $G:\ \eta\in X^+_*(H)$ iff
$\langle\eta,\check\alpha_i\rangle\geq0$ for any $\check\alpha_i\in \Delta$.
For example, if $H=G=\bH$, then $X^+_*(H)=Y$. If $\phi\in\CM(\bX,\alpha)$,
we may view $\phi$ as a reduction of the trivial $\bG$-torsor to $\bB$.
Let $\phi^\bH$ be the corresponding induced $\bH$-torsor. Then its isomorphism
class equals $-\alpha$.

Let us view $\psi$ as a reduction of the trivial $\bG$-torsor to $\bP_J$.
Let $\bL_J$ be the Levi quotient of $\bP_J$, and let $\psi^{\bL_J}$ be the
corresponding induced $\bL_J$-torsor.
Let $\varphi$ be the {\em Harder-Narasimhan flag} of $\psi^{\bL_J}$. We may view
it as a reduction of $\psi$ to a parabolic subgroup $\bP_K,\ K\subset J$.
By definition, the isomorphism class $\eta$ of $\varphi^{\bL_K}$
(as an element of $Y$) has the following properties:

a) $\varpi_*\eta=-\beta$;

b) $\langle\eta-\varpi_*\eta,\check\alpha_j\rangle>0$ for $j\in J-K$;

c) $\langle\eta-\varpi_*\eta,\check\alpha_k\rangle=0$ for $k\in K$.

In particular, if $\bL'_K$ stands for the quotient of $\bL_K$ by its center,
then the induced torsor $\varphi^{\bL'_K}$ is trivial. Choosing its trivial
reduction to the positive Borel subgroup of $\bL'_K$ we obtain a reduction
$\phi$ of $\psi$ to $\bB$. Thus $\phi$ is a map from $C$ to $\bX$ of
degree $\alpha=-\eta$. We see that $\alpha$ is a special lift of $\beta$,
and $\phi\in\CM_b(\bX,\alpha)$ has the desired properties.
\end{proof}

\end{document}